\hsize=11.4cm
\vsize=18cm
\hoffset=3cm
\voffset=3cm

\parindent 20pt
\tolerance=10000

\centerline{\bf Compressions and Pinchings}
\bigskip
\centerline{Jean-Christophe Bourin}
\medskip
\centerline{Les Coteaux, rue Henri Durel 78510 Triel, France}
\smallskip
\centerline{E-mail: bourinjc@club-internet.fr}

\bigskip
\bigskip
\noindent
{\bf Abstract} \quad {\it  There exist operators $A$ such that : for any sequence of
contractions $\{ A_n\}$, there is a total sequence of mutually orthogonal projections $\{ E_n\}$ 
such that $\Sigma E_nAE_n=\bigoplus A_n$. }


\bigskip
\medskip
\noindent
{\bf Introduction}
\bigskip

 By an operator, we mean an element in the algebra ${\rm L}({\cal H})$ of all bounded linear
operators acting on the usual (complex, separable, infinite dimensional) Hilbert space 
${\cal H}$. We denote by the same letter a projection and the corresponding subspace.
If $F$ is a projection and $A$ is an operator, we denote by $A_F$ the compression of $A$ by $F$,
that is the restriction of $FAF$ to the subspace $F$. Given a
total sequence of nonzero mutually orthogonal projections $\{E_n\}$, we consider the pinching
$$
{\cal P}(A)=\sum_{n=1}^{\infty}E_nAE_n=\bigoplus_{n=1}^{\infty}A_{E_n}.
$$
If $\{A_n\}$ is a sequence of operators acting on separable Hilbert spaces with $A_n$ 
unitarily equivalent to $A_{E_n}$ for all $n$, we also naturally write 
${\cal P}(A)=\bigoplus_{n=1}^{\infty}A_n$. The main result of this paper
can then be stated as :
\bigskip
\noindent
\qquad\qquad {\it Let $\{A_n\}_{n=1}^{\infty}$ be a sequence of operators
acting on separable Hilbert spaces. Assume that $\sup_n\Vert A_n \Vert<1$.
Then, we have a pinching
$$
{\cal P}(A)=\bigoplus_{n=1}^{\infty}A_n
$$
for any operator $A$ whose essential numerical range contains the unit disc. }
\bigskip
\noindent
This result is proved in the second section of the paper. We have included a first section
concerning some well-known properties of the essential numerical range.

\bigskip\medskip\noindent
{\bf 1. Properties of the essential numerical range }

\bigskip
We denote by $\langle\cdot,\cdot\rangle$ the inner product (linear in the second variable),
by ${\rm co}S$ the convex hull of a subset $S$ of the complex plane ${\rm \bf C}$.
$W(A)=\{ \langle h,Ah\rangle\ |\quad\Vert h\Vert=1 \}$ is the numerical range of the 
operator $A$ and  $\overline{W}(A)$ is the closure of $W(A)$.
The celebrated Hausdorff-Toeplitz theorem (cf [6] chapter 1) states that $W(A)$ is convex. 
A corollary is Parker' s theorem ([6], p.20) : Given a $n$ by $n$ matrix $A$, there is a 
matrix $B$ unitarily equivalent to $A$ and with all its diagonal elements equal to
 ${\rm Tr}A/n$.
 
Here are three equivalent definitions of the essential numerical range of $A$,
denoted by $W_e(A)$ :

\bigskip\noindent
(1)\qquad $W_e(A)=\cap\overline{W}(A+K)$
where the intersection runs over the compact operators $K$ 
\bigskip
\noindent
(2)\qquad  Let $\{E_n\}$ be any sequence of finite rank projections converging strongly to the 
identity and denote by $B_n$ the compression of $A$ to the subspace $E_n^{\perp}$. Then
$W_e(A)=\cap_{n\ge1}\overline{W}(B_n)$ 
\bigskip\noindent
(3)\qquad $W_e(A)=$
$$
\{\ \lambda\ |\quad {\rm there\ is\ an\ orthonormal\ system\ \{ e_n \}_{n=1}^\infty \
with\ \lim\langle e_n, Ae_n\rangle =\lambda} \}.
$$
It follows that $W_e(A)$ is a compact convex set containing the essential spectrum of $A$,
$Sp_e(A)$. The equivalence between these definitions has been known since the early seventies
if not sooner (see for instance [1]). The very first definition of $W_e(A)=$ is (1); however (3)
is also a natural notion and easily entails convexity and compactness of the essential numerical
range. We mention the following result of Chui-Smith-Smith-Ward [4] :
 
\bigskip\noindent
{\bf Proposition 1.} {\it Every operator $A$ admits some compact perturbation $A+K$
for which $W_e(A)=\overline{W}(A+K)$. }  

\bigskip
Another characterization of the essential numerical range of $A$ is 
$$
W_e(A)
=\{\ \lambda\ |\quad {\rm there\ is\ a\ basis\ \{ e_n \}_{n=1}^\infty \
 with\ \lim\langle e_n, Ae_n\rangle =\lambda} \}.
$$
Let us check the equivalence between our definition (3) with orthonormal system and the
above identity which seems to be due to Q. F. Stout [11].
Let $\{ x_n \}_{n=1}^\infty$ be an orthonormal system such that
$\lim_{n\to\infty}\langle x_n, Ax_n\rangle=\lambda$. 
If ${\rm span}\{ x_n \}_{n=1}^\infty$ is of finite codimension $p$ we immediately get a basis
$e_1,\dots,e_p;e_{p+1}=x_1,\dots,e_{p+n}=x_n,\dots$ such that 
$\lim_{n\to\infty}\langle e_n, Ae_n\rangle=\lambda$. 
If ${\rm span}\{ x_n \}_{n=1}^\infty$ is of infinite codimension, we may complete this system
with $\{y_n\}_{n=1}^\infty$ in order to obtain a basis. Let $P_j$ be the subspace spanned by
$y_j$ and $\{x_n|2^{j-1}\le n< 2^j\}$. By Parker's theorem, there is a basis of $P_j$, 
say $\{e_l^j\}_{l\in\Lambda_j}$, with
$$
\langle e_l^j,Ae_l^j\rangle={1\over {\rm dim}P_j}{\rm Tr}AP_j.
$$
Since
$$
{1\over {\rm dim}P_j}{\rm Tr}AP_j\rightarrow\lambda\quad{\rm as}\quad j\rightarrow\infty,
$$
we may index $\{e_l^j\}_{j\in{\rm N};l\in\Lambda_j}$ in order to obtain a basis 
$\{ f_n \}_{n=1}^\infty$ such that $\lim_{n\to\infty}\langle f_n, Af_n\rangle=\lambda$. 

\bigskip
The essential numerical range appears closely related to the diagonal set of $A$ which we define by
$$
\Delta(A)
=\{\ \lambda\ |\quad {\rm there\ is\ a\ basis\ \{ e_n \}_{n=1}^\infty \
 with\ \langle e_n, Ae_n\rangle =\lambda} \}.
$$
The next result is a straightforward consequence of a lemma of Peng Fan [5]. A real operator
 means an operator acting on a real Hilbert space and int$X$ denotes the interior of $X\subset
{\rm\bf C}$.

\bigskip
\noindent
{\bf Proposition 2.} {\it  Let $A$ be an operator. Then 
${\rm int}\,W_e(A)\subset\Delta(A)\subset W_e(A)$. Consequently, an open set ${\cal U}$ is
contained in $\Delta(A)$ if and only if there is a basis $\{e_n\}_{n=1}^{\infty}$ such that
${\cal U}\subset{\rm co}\{\langle e_k,Ae_k\rangle|k\ge n\}$ for all $n$.
Finally, the diagonal set of a real operator is  symetric about the real axis.
(For $A$  self-adjoint, the result holds with ${\rm int}$ denoting
the interior of subsets of ${\rm\bf R}$.)  }

\bigskip\noindent
Curiously enough, it seems difficult to answer the following questions :
 Is the diagonal set always a (possibly vacuous) convex set ? Is there an operator of the form
self-adjoint + compact with  a disconnected diagonal set ?

An elementary but very important property of $W(\cdot)$ is the so named projection property :
${\rm Re}W(A)=W({\rm Re}A)$ (see [6] p. 9), where Re stands for real part. $W_e(\cdot)$ has also this property. 
 This result and the Hausdorff-Toeplitz Theorem are the keys to prove the following fact :

\bigskip\noindent
{\bf Proposition 3.} {\it Let $A$ be an operator. 
\smallskip\noindent
{\rm (1)} If $W_e(A)\subset W(A)$ then $W(A)$ is closed.
\smallskip\noindent
{\rm (2)} There exist normal, finite rank operators $R$ of
arbitrarily small norm such that $W(A+R)$ is closed. }

\bigskip\noindent
{\bf Proof.} Assertion (1) is due to J. S. Lancaster [8]. We prove the second assertion and 
implicitly prove Lancaster's result.
We may find an orthonormal system $\{f_n\}$ such that the closure of the sequence
$\{\langle f_n,Af_n\rangle\}$ contains the boundary $\partial W_e(A)$. Fix $\varepsilon>0$. 
It is possible to find an integer p and scalars $z_j$, $1<j<p$, with $|z_j|<\varepsilon$,
such that :
$$
{\rm co}\{ \langle f_j,Af_j\rangle+z_j\ |1<j<p\} \supset \partial W_e(A). 
$$
Thus, the finite rank operator $R=\sum_{1<j<p}z_j f_j\otimes f_j$ has the property that
${\rm W}(A+R)$ contains $W_e(A)$. 

 We need this  operator $R$. Indeed, setting $X=A+R$, we also have
${\rm W}(X)\supset W_e(X)$. We then claim that ${\rm W}(X)$ is closed (this claim
implies assertion (1)). By the contrary, there would exist 
$z\in\partial\overline{W}(X)\setminus W_e(X)$. Furthermore, since
$\overline{W}(X)$ is the convex hull of its extreme points, we could assume that such 
a z is an extreme of $\overline{W}(X)$. By suitable rotation and 
translation, we could assume that $z=0$ and that the imaginary axis is a line of support of
$\overline{W}(X)$. The projection property for $W(\cdot)$ would imply that
$W({\rm Re}X)=(x,0[$ for a certain negative number $x$, so that
$0\in W_e({\rm Re}X)$. Thus we would deduce from the projection property for $W_e(\cdot)$ that
$0\in W_e(X)$ : a contradiction.  \quad ${\bf \diamondsuit}$

\bigskip\noindent
 The perturbation $R$ in Proposition 4 can be taken real if $A$ is real. We mention that the set
of operators with nonclosed numerical ranges is not dense in  ${\rm L}({\cal H})$. 
Proposition 3 improves the 
following result of I.D. Berg and B. Sims [3] : operators which attain their numerical radius 
are norm dense in ${\rm L}({\cal H})$. A motivation for Berg and Sims was the following fact 
 :  given an arbitrary operator $A$, a small rank one perturbation
of $A$ yields an operator which attains its norm. Indeed, the polar decompositon allows us to assume that
$A$ is positive, an easy case when reasoning as in the proof of Proposition 3.  

Let us say that a convex set in ${\bf C}$ is relatively open if either it is a single point, 
an open segment or an usual open set. Using similar methods as in the previous proof, or 
applying Propositions 2 and 3, we obtain 

\bigskip\noindent
{\bf Proposition 4.} {\it For an operator $A$ the following assertions are equivalent
\smallskip\noindent
{\rm (a)}\quad $W(A)$ is relatively open.
\smallskip\noindent
{\rm (b)}\quad $\Delta(A)=W(A)$. }

\bigskip
>From the previous results we may derive some information about $W(\cdot)$, $W_e(\cdot)$ and
$\Delta(\cdot)$ for various classes of operators :  
\smallskip\noindent
{\rm (a)}\quad Let $S$ be either the unilateral or bilateral Shift, then
 $\Delta(S)=W(S)$ is the open unit disc. More generally Stout showed [10] that weighted
periodic shifts $S$ have open numerical ranges; therefore $\Delta(S)=W(S)$. 
\smallskip\noindent
{\rm (b)}\quad There exist a number of Toeplitz operators with open numerical range. See the papers
by E. M. Klein [7] and by J. K. Thukral [11].
\smallskip\noindent
{\rm (c)}\quad  Let $X$ be an operator lying in a $C^*$-subalgebra of $L({\cal H})$ with no
finite dimensional projections. Then for any real $\theta$,
$\overline{W}({\rm Re}e^{{\rm i}\theta}X)= W_e({\rm Re}e^{{\rm i}\theta}X)$. From the projection
property for $W(\cdot)$ and  $W_e(\cdot)$ we infer that $W_e(X)=\overline{W}(X)$. 
\smallskip\noindent
{\rm (d)}\quad Let $X$ be an essentially normal operator : $X^*X-XX^*$ is compact. It is known that
$W_e(X)={\rm co}Sp_e(X)$. Indeed, for such an operator the essential norm equals the essential 
spectral radius : $\Vert X\Vert_e=\rho_e(X)$. Denoting by $w_e(X)$ the essential numerical radius
of $X$ we deduce that $\Vert X\Vert_e=w_e(X)=\rho_e(X)$. Note that $e^{{\rm i}\theta}X+\mu I=Y$
is also an essentially normal operator for any $\theta\in {\rm \bf R}$ and $\mu\in {\rm \bf C}$.
Let $z$ be an extremal point of $W_e(X)$. With suitable $\theta$ and $\mu$ we have
$e^{{\rm i}\theta}z+\mu=w_e(Y)=\max\{|y|,\ y\in W_e(Y)\}$, the maximum being attained at the
single point $e^{{\rm i}\theta}z+\mu$. Since ${\rm co}Sp_e(Y)\subset W_e(Y)$ and
$\rho_e(Y)=w_e(Y)$, this implies that  $e^{{\rm i}\theta}z+\mu \in Sp_e(Y)$. Hence
$z\in Sp_e(Y)$, so that $W_e(X)={\rm co}Sp_e(X)$.

\bigskip
\medskip
\noindent
{\bf 2. The pinching theorem}

\bigskip
Recall that one way to define the essential numerical range of an operator $A$ is : $W_e(A)=$
$$
\{\ \lambda\ |\quad {\rm there\ is\ an\ orthonormal\ system\ \{ e_n \}_{n=1}^\infty \
with\ \lim\langle e_n, Ae_n\rangle =\lambda} \}.
$$
It is then easy to check that $W_e(A)$ is a compact convex set.
Moreover $W_e(A)$ contains the open unit disc ${\cal D}$ if and only if
there is a basis $\{e_n\}_{n=1}^{\infty}$ such that ${\rm co}\{\langle e_k,Ae_k\rangle
\ |\ k>n\}\supset{\cal D}$ for all $n$. 

\bigskip
\noindent
{\bf Theorem 1.} {\it Let $A$ be an operator with $W_e(A)\supset{\cal D}$
and let $\{A_n\}_{n=1}^{\infty}$ be a sequence of operators such that
$\sup_n\Vert A_n\Vert<1$. Then, we have a pinching
$$
{\cal P}(A)=\bigoplus_{n=1}^{\infty}A_n.
$$
$($If $A$ and $\{A_n\}_{n=1}^{\infty}$ are real, then we may take a real pinching$)$.
}  

\bigskip
\noindent
{\bf Proof.} It suffices to solve the following problem
\bigskip
\noindent
[P]\quad Let $A$ be an operator, with
$\Vert A\Vert\le\gamma$ and $W_e(A)\supset{\cal D}$, let $h$ be a norm one vector 
and $X$ a strict contraction, $\Vert X\Vert<\rho<1$. Find
a projection $E$, and a constant $\varepsilon>0$ only depending on
$\gamma$ and $\rho$ such that :
\smallskip
\ (i)\quad $\dim E=\infty$ and $A_E=X$
\smallskip
(ii)\quad $\dim E^{\perp}=\infty$, $W_e(A_{E^{\perp}})\supset{\cal D}$ and $\Vert Eh\Vert\ge
\varepsilon$.
\bigskip
\noindent
Let us explain why it is sufficient to solve [P]. Take $\gamma=\Vert A\Vert$ and fix a dense
sequence $\{ h_n\}$ in the unit sphere of ${\cal H}$. We claim that (i) and (ii) ensure that
there exists a sequence of mutually orthogonal projections $\{E_j\}$ such that, setting
$F_n=\sum_{j\le n}E_j$, we have for all integers $n$ :
\smallskip\noindent
\quad \ ($*$) $A_n=A_{E_n}$ and $W_e(A_{F_n^{\perp}})\supset{\cal D}$ \
(so $\dim F_n^{\perp}=\infty$)
\smallskip\noindent
\quad ($**$) $\Vert F_n h_n\Vert \ge \varepsilon$.
\smallskip\noindent
 This is true for $n=1$ by (i). Suppose this holds for an $N\ge1$. Let
$\nu(N)\ge N+1$ be the first integer for which $F_N^{\perp} h_{\nu(N)} \neq 0$. Note that
$\Vert A_{F_N^{\perp}}\Vert \le \gamma$. We apply (i) and (ii) to $A_{F_N^{\perp}}$,
$A_{N+1}$ and $F_N^{\perp} h_{\nu(N)} / \Vert F_N^{\perp} h_{\nu(N)}\Vert$ in place of $A$, $X$ and $h$. 
We then deduce that $(*)$ and $(**)$ are still valid for $N+1$. Therefore
$(*)$ and $(**)$ hold for all $n$. Denseness of $\{ h_n\}$ and $(**)$ show that $F_n$ strongly
increases to the identity $I$ so that $\sum_{j=1}^{\infty}E_j=I$ as required.
 
\bigskip
We first solve problem [P] restricted to condition (i), consisting in 
representing $A$ as a dilation of $X$. 
Next, we solve problem [P] completely.
\bigskip
\noindent
\centerline{ {\it 1. Preliminaries} }
\medskip
\noindent
We shall use a sequence $\{V_k\}_{k\ge1}$ of orthogonal matrices acting on spaces
of dimensions $2^k$. This sequence is built up by induction :
$$
V_1={1\over\sqrt2}\pmatrix{1&1\cr-1&1}\quad{\rm then}\quad
V_k={1\over\sqrt2}\pmatrix{V_{k-1}&V_{k-1}\cr-V_{k-1}&V_{k-1}}\quad{\rm for}\ k\ge2.
$$
Given a Hilbert space ${\cal G}$ and a decomposition
$$
{\cal G}=\bigoplus_{j=1}^{2^k}{\cal H}_j\quad {\rm with}\ {\cal H}_1=\dots={\cal H}_{2^k}
={\cal H},
$$
we may consider the unitary
(orthogonal) operator on ${\cal G}$ : $W_k=V_k\bigotimes I$,
where $I$ denotes the identity on ${\cal H}$, 

Now, let $B:{\cal G}\rightarrow{\cal G}$ be an operator which, relatively to the above
decomposition of ${\cal G}$, is written with a block diagonal matrix
$$
B=\pmatrix{B_1&\ &\ \cr \ &\ddots&\ \cr\ &\ &B_{2^k}}.
$$
We observe that the block matrix representation of $W_kBW_k^*$ has its diagonal
entries all equal to
$$
{1\over2^k}\left(B_1+\dots B_{2^k}\right).
$$
So, the orthogonal operators $W_k$ allow us to pass from a block diagonal matrix representation
to a block matrix representation in which
the diagonal entries are all equal.
\bigskip
\noindent
\centerline{ {\it 2. Solution of problem} [P]-(i) }
\medskip
\noindent
The contraction $Y=(1/\Vert X\Vert)X$ can be dilated in a unitary 
$$
U=\pmatrix{Y&-(I-YY^*)^{1/2}\cr (I-Y^*Y)^{1/2}&Y^*}
$$
thus $X$ can be dilated in a normal operator $N=\Vert X\Vert U$ with
$\Vert N\Vert<\rho$.
This permits to restrict  to the case when $X$ is a normal contraction,
$\Vert X\Vert<\rho<1$. Thus we set the following problem :
\bigskip
\noindent
[Q]\quad Let $X$ be a normal contraction, $\Vert X\Vert<\rho<1$. Find 
a projection $E$, $\dim E=\infty$, such that $A_E=X$

\bigskip
\noindent
We remark with the Berg-Weyl-von Neumann theorem [2] that a normal contraction $X$,
 $\Vert X\Vert<\rho<1$, can be written
$$
X=D+K \eqno (1)
$$
where $D$ is normal diagonalizable, $\Vert D\Vert=\Vert X\Vert<\rho$, and $K$ is compact with an
arbitrarily small norm. Let $K={\rm Re}K+{\rm iIm}K$ be the cartesian decomposition of $K$. 
We can manage to have an integer $l$, a real $\alpha$
and a real $\beta$ so that decomposition (1) satisfies :

\noindent
a) the operators $\alpha D$, $\beta{\rm Re}K$, $\beta{\rm Im}K$ are majorized in norm by $\rho$,

\noindent
b) there are positive integers $m$, $n$ with $2^l=m+2n$ and
$$
X={1\over2^l}(m\alpha D+n\beta{\rm Re}K+n\beta{\rm iIm}K). \eqno (2)
$$
More precisely we can take any $l$ such that $[2^l/(2^l-2)].\Vert X\Vert<\rho$. Next, assuming
$\Vert K\Vert<\rho/2^l$, we can take $m=2^l-2$, $n=1$, $\alpha=2^l/(2^l-2)$ and $\beta=2^l$.

Let then $T$ be the diagonal normal operator acting on the space
$$
{\cal G}=\bigoplus_{j=1}^{2^l}{\cal H}_j\quad {\rm with}\ {\cal H}_1=\dots={\cal H}_{2^l}
={\cal H},
$$
and defined by
$$
T=\left( \bigoplus_{j=1}^m D_j\right)\bigoplus
\left( \bigoplus_{j=m+1}^{m+n} R_j\right)\bigoplus
\left( \bigoplus_{j=m+n+1}^{2^l} S_j\right)
$$
where $D_j=\alpha D$, $S_j=\beta{\rm Re}K$ and $S_j=\beta{\rm iIm}K$.
\smallskip
\noindent
We note that $\Vert T\Vert<\rho<1$ and that the operator $W_lTW_l^*$, represented in the
preceding decomposition of ${\cal G}$, has its diagonal entries all equal to
$X$ by (2). Thus to solve [Q] it suffices to solve the following problem
\smallskip
\noindent
[R]\quad Given a diagonal normal operator $T$, $\Vert T\Vert<\rho<1$,
find a projection $E$, $\dim E=\infty$, such that $A_E=T$.
\smallskip
\noindent
{\it Solution of }[R]. Let $\{ \lambda_n(T)\}_{n\ge1}$ be the eigenvalues of $T$ repeated
according to their multiplicities. Since $|\lambda_n(T)|<1$ for all $n$ and that $W_e(A)
\supset{\cal D}$, we have a norm one vector $e_1$ such that $\langle e_1,Ae_1\rangle=
\lambda_1(T)$. Let $F_1=[{\rm span}\{e_1,Ae_1,A^*e_1\}]^{\perp}$. As $F_1$ is of finite 
codimension, $W_e(A_{F_1})\supset{\cal D}$. So, there exists a norm one vector $e_2\in F_1$
such that $\langle e_2,Ae_2\rangle=\lambda_2(T)$. Next, we set
$F_2=[{\rm span}\{e_1,Ae_1,A^*e_1,e_2,Ae_2,A^*e_2\}]^{\perp}$, $\dots$. If we go on like this, 
we exhibit an orthonormal system $\{e_n\}_{n\ge1}$ such that, setting 
$E={\rm span}\{e_n\}_{n\ge1}$, we have $A_E=T$. \quad ${\bf \diamondsuit}$ 
 
\bigskip
\noindent
\centerline{ {\it 3. Solution of problem} [P], (i) {\it and} (ii) }
\medskip
\noindent
We take an arbitrary norm one vector $h$. We can show, using the same reasoning as that applied
to solve [R], that we have an orthonormal system $\{f_n\}_{n\ge0}$,
with $f_0=h$, such that :

\noindent
a)\quad $\langle f_{2j}, Af_{2j}\rangle=0$ for all $j\ge1$

\noindent
b)\quad $\{ \langle f_{2j+1}, Af_{2j+1}\rangle \}_{j\ge0}$ is a dense sequence in ${\cal D}$ 

\noindent
c)\quad If $F={\rm span}\{f_j\}_{j\ge0}$, then $A_F$ is the normal operator
$$
\sum_{j\ge0}\langle f_j, Af_j\rangle f_j\otimes f_j.
$$
Setting $F_0={\rm span}\{f_{2j}\}_{j\ge0}$ and $F'_0={\rm span}\{f_{2j+1}\}_{j\ge0}$, we then have :

\noindent 
a)\quad Relatively to the decomposition $F=F_0\bigoplus F'_0$, $A_F$ can be written
$$
A_F=\pmatrix{A_{F_0}&0\cr 0&A_{F'_0}}.
$$

\noindent 
b)\quad $W_e(A_{F'_0})\supset{\cal D}$ and $h\in F_0$.

\noindent
We can then write a decomposition of  $F'_0$, $F'_0=\bigoplus_{j=1}^{\infty}F_j$
where for each index $j$, $F_j$ commutes with $A_F$ and $W_e(A_{F_j})\supset{\cal D}$;
so that the decomposition $F=\bigoplus_{j=0}^{\infty}F_j$ yields a representation
of $A_F$ as a block diagonal matrix, 
$$
A_F=\bigoplus_{j=0}^{\infty}A_{F_j}.
$$
Since $W_e(A_{F_j})\supset{\cal D}$ when $j\ge1$, the same reasoning as that used in the 
solution of [R] entails that for any sequence  
$\{ X_j \}_{j\ge0}$ of strict contractions we have decompositions
$(\dagger)$ $F_j=G_j\bigoplus G'_j$ allowing us to write $(j\ge1)$
$$  
A_{F_j}=\pmatrix{X_j&*\cr *&*}.
$$
Since $\Vert X\Vert<\rho<1$, we can find an integer $l$ only depending on $\rho$
and $\gamma$, as well as strict contractions $X_1,\dots,X_{2^l}$, such that
$$
X={1\over2^l}\left( A_{F_0}+\sum_{j=1}^{2^l-1}X_j \right). \eqno (3)
$$
We come back to decompositions $(\dagger)$ and we set
$$
G=F_0\bigoplus \left( \bigoplus_{j=1}^{2^l-1}G_j \right). 
$$
Relatively to this decomposition,
$$
A_G=\pmatrix{A_{F_0}&\ &\ &\ \cr
\ &X_1&\ &\ \cr
\ &\ &\ddots&\ \cr
\ &\ &\ &X_{2^l-1}}.
$$ 
Then we deduce from (3) that the block matrix $W_lA_GW_l^*$ has its diagonal entries 
all equal to $X$.

Summary :
 $h\in G$ and there exists a decomposition $G=\bigoplus_{j=1}^{2^l}E_j $
such that $A_{E_j}=X$ for each $j$. Thus we have an integer $j_0$ such that, setting 
$E_{j_0}=E$, we have
$$
A_E=X \quad {\rm and} \quad \Vert Eh\Vert\ge {1\over\sqrt{2^l}}.                                             
$$ 
The proof is finished. ${\bf \diamondsuit}{\bf \diamondsuit}$

\bigskip
\noindent
{\bf Corollary 1.} {\it Let $A$ be an operator with $W_e(A)\supset{\cal D}$. For any strict
contraction $X$, there is an isometry $V$ such that $X=V^*AV$. }

\bigskip
\noindent
{\bf Corollary 2.} {\it Let $A$ be an operator with $W_e(A)\supset{\cal D}$. For any 
contraction $X$, there is a sequence  $\{U_n\}$ of unitary operators such that
$U_n^*AU_n\rightarrow X$ in the weak operator topology. }

\bigskip
We use the strict inclusion notation $X\subset\subset Y$ for subsets $X$, $Y$ of ${\rm \bf C}$
to mean that there is an $\varepsilon>0$ such that 
$\{ x+z\ \ |\ x\in X,\ |z|<\varepsilon \}\subset Y$.

\bigskip
\noindent
{\bf Theorem 2.} {\it Let $A$ be an operator and let $\{ A_n \}_{n=1}^{\infty}$ be 
a sequence of normal operators. If $\cup_{n=1}^{\infty}W(A_n)\subset\subset W_e(A)$ 
then we have a pinching
$$
{\cal P}(A)=\bigoplus_{n=1}^{\infty}A_n.
$$
(For self-adjoint operators, this result holds with the strict inclusion of ${\rm\bf R}$.)
}

\bigskip
\noindent
{\bf Sketch of proof.} Let $N$ be a normal operator with $W(N)\subset\subset W_e(A)$. If
$N$ is diagonalizable, reasonning as in the proof of Theorem 1, we deduce that $N$ can
be realized as a compression of $A$. If $N$ is not diagonalizable we may assume that
$0\in W_e(A)$. Thanks to the Berg-Weyl-von Neumann Theorem and still reasonning as in the proof
of Theorem 1 we again deduce that $N$ is a compression of $A$. Finally, the strict containment
assumption allows us to get the wanted pinching. 

\bigskip
To finish this section, we mention that we can not drop the assumption that the strict
contractions $A_n$ of Theorem 1 are uniformly bounded in norm by a real $<1$.
This observation is equivalent to the fact that we can not delete the strict containment
assumption in Theorem 2 :

 Let $P$ be a halving projection ($\dim P=\dim P^{\perp}=\infty$),
 so $W_e(P)=[0,1]$. Then the sequence $\{ 1-1/n^2\}_{n\ge1}$
can not be realized as the entries of the main diagonal of a matrix representation of $P$.
 To check that, we note that the positive operator $I-P$ would be in the trace-class :
a contradiction. (Recall that a positive operator with a summable diagonal is trace class.)

\bigskip\bigskip
\centerline{\bf References}
\bigskip\noindent
[1] J.H. Anderson and J.G. Stampfli, {\it Commutators and compression}, Israel J. Math {\bf 10}
(1971), 433-441

\noindent
[2] I.D. Berg, {\it A extension of the Weyl-von Neumann theorem to normal operators} Trans. Amer.
Math. Soc. {\bf 160} (1971), 365-371

\noindent
[3] I.D. Berg and B. Sims, {\it Denseness of operators which attain their numerical ranges}, 
J. Austral. Math. Soc. {\bf 36}, serie A (1984), 130-133

\noindent
[4] C.K. Chui, P.W. Smith, R.R. Smith, J.D. Ward, {\it L-ideals and numerical range preservation},
Illinois J. Math. {\bf 21}, 2 (1977), 365-373.

\noindent
[5] P. Fan, {\it On the diagonal of an operator}, Trans. Amer.Math. Soc. {\bf 283} (1984), 
239-251

\noindent
[6] R.A. Horn and C.R. Johnson, {\it Topics in matrix analysis}, Cambridge University Press, 1990

\noindent
[7] E.M. Klein, {\it The numerical range of a Toeplitz operator}, Proc. Amer. Math. Soc. {\bf35}
(1972), 101-103

\noindent
[8] J.S. Lancaster, {\it The boundary of the numerical range}, Proc. Amer. Math. Soc. {\bf49}, 2
(1975), 393-398

\noindent
[9] Q.F. Stout, {\it Shur products of operators and the essential numerical range},
Trans. Amer. Math. Soc {\bf 264}, 1 (1981), 39-47

\noindent
[10] Q.F. Stout, {\it The numerical range of  a weighted shift}, Proc. Amer. Math. Soc.
{\bf 88}, 3 (1983), 495-502

\noindent
[11] J.K. Thukral, {\it The numerical range of a toeplitz operator with harmonic symbol},
J. Operator theory {\bf 34} (1995), 213-216

\bye